\documentclass[reqno]{amsart}
\usepackage{amsthm, amssymb, amsfonts}
\usepackage{lmodern}
\usepackage{color} 
\usepackage{hyperref}
\usepackage[T5]{fontenc}
\usepackage{amscd,amssymb}
\usepackage[v2,cmtip]{xy}
\usepackage{mathrsfs}
\theoremstyle{plain}
\newtheorem{thm}{Theorem}[section]

\newtheorem{con}[thm]{Conjecture}
\newtheorem{corl}[thm]{Corollary}
\theoremstyle{definition}
\newtheorem{defn}[thm]{Definition}
\newtheorem{rem}[thm]{Remark}

\theoremstyle{plain}

\theoremstyle{definition}

\begin{document}  

\title[On the hit problem for the Steenrod algebra] {On the hit problem for the Steenrod algebra \\in some generic degrees and applications} 

\author[N. K. Tin]{Nguyen Khac Tin} 

\address{Faculty of Applied Sciences, 
Ho Chi Minh City University of Technology and Education \\
01 Vo Van Ngan, Thu Duc, Ho Chi Minh city, Viet Nam.}

\email{tinnk@hcmute.edu.vn}

\maketitle

\begin{abstract}
Let $\mathcal P_{n}:=H^{*}((\mathbb{R}P^{\infty})^{n}) \cong \mathbb F_2[x_{1},x_{2},\ldots,x_{n}]$ be the polynomial algebra over the prime field of two elements, $\mathbb F_2.$ We investigate the Peterson hit problem for the polynomial algebra $\mathcal P_{n},$ viewed as a graded left module over the mod-$2$ Steenrod algebra, $\mathcal{A}.$ For $n>4,$ this problem is still unsolved, even in the case of $n=5$ with the help of computers. 

The purpose of this paper is to continue our study of the hit problem by developing a result in \cite{ph31} for $\mathcal P_n$ in the generic degree $r(2^s-1)+m.2^s$ where $r=n=5,\ m=13,$ and $s$ is an arbitrary non-negative integer. Note that for $s=0,$ and $s=1,$ this problem has been studied by Phuc \cite{ph20ta}, and \cite{ph31}, respectively. 

As an application of these results, we get the dimension result for the polynomial algebra in the generic degree $d=(n-1).(2^{n+u-1}-1)+\ell.2^{n+u-1}$ where $u$ is an arbitrary non-negative integer, $\ell \in \{23, 67 \},$ and $n=6.$

One of the major applications of hit problem is in surveying a homomorphism introduced by Singer, which is a homomorphism
$$Tr_n :\text{Tor}^{\mathcal A}_{n, n+d} (\mathbb F_2,\mathbb F_2) \longrightarrow  (\mathbb{F}_{2}{\otimes}_{\mathcal{A}}\mathcal P_{n})_d^{GL(n; \mathbb F_2)}$$
from the homology of the Steenrod algebra to the subspace of $(\mathbb{F}_{2}{\otimes}_{\mathcal{A}}\mathcal P_{n})_d$ consisting of all the $GL(n; \mathbb F_2)$-invariant classes. It is a useful tool in describing the homology groups of the Steenrod algebra, $\text{Tor}^{\mathcal A}_{n, n+d}(\mathbb F_2,\mathbb F_2).$
The behavior of the fifth Singer algebraic transfer in degree $5(2^s-1)+13.2^s$ was also discussed at the end of this paper. 
\end{abstract}

\medskip
\section{Introduction}\label{s1} 
\setcounter{equation}{0}

Throughout the paper, we denote a prime field with two elements by $\mathbb F_2.$ Let $\mathbb{R}\mathcal P^{\infty}$ be the infinite dimensional real projective spaces. Then, $H^{*}(\mathbb{R}\mathcal P^{\infty}) \cong \mathbb F_2[x_{1}],$ and therefore, the mod-$2$ cohomology algebra of the direct product of $n$ copies of $\mathbb{R}\mathcal P^{\infty}$ is isomorphic to the graded polynomial algebra 
$\mathbb F_2[x_{1},x_{2},\ldots,x_{n}],$ reviewed as an unstable $\mathcal{A}$-module on $n$ generators $x_1, x_2, \ldots, x_n,$ each of degree one.

The action of $\mathcal{A}$ on $\mathcal P_n$ is determined by the formula 
\[
Sq^{i}(x_{j})=\left\{ \begin{array}{ll}
x_{j}, & i=0,\\
x_{j}^{2}, & i=1,\\
0, & i>1,
\end{array}\right.
\]
and the Cartan formula $Sq^{k}(uv)=\sum_{i=0}^{k}Sq^{i}(u)Sq^{k-i}(v),$ where $u, v \in \mathcal P_n$ (see Steenrod and Epstein~\cite{s-e62}).  

One of central problems of Algebraic Topology is to determine a minimal set of generators for the ring of invariants $(\mathcal P_n)^{G_n},$ where $G_n \subset GL(n; \mathbb F_2),$ the general linear group of invertible matrices. This ring is stable under the action of $\mathcal{A}.$ The problem is so-called hit problem for the polynomial algebra. If we consider $\mathbb{F}_{2}$ as a trivial $\mathcal{A}$-module, then the hit problem is equivalent to the problem of finding a basis for the $\mathbb F_2$-graded vector space 
\begin{align}\label{ct1}
\left\{(\mathbb{F}_{2}{\otimes}_{\mathcal{A}}(\mathcal P_{n})^{G_n})_d\right\}_{d \geqslant 0}.
\end{align}

The structure of this space has first been studied by Peterson \cite{pe87}, Wood \cite{wo89}, Singer \cite {si89}, Priddy \cite{pr90},  who show its relationship to several classical problems in cobordism theory, modular respresentation theory, Adams spectral sequence for the stable homotopy of spheres, stable homotopy type of the classifying space of finite groups. Then, this problem was investigated by Hung-Peterson~\cite{h-p95}, Kameko~\cite{ka98}, Janfada~\cite{ja07, ja09}, Nam~\cite{na04}, Repka-Selick \cite{r-s98}, Silverman~\cite{sil95}, Wood~\cite{wo89}, Sum~\cite{su10, su15}, Sum-Tin~\cite{s-t21, t-s16}, the present writer~\cite{ti21pja}, \cite{ti21aejm} and others. 

For $G_n=\Sigma_n$ the symmetric group on $n$ letters, the space (\ref{ct1}) is known by Janfada-Wood~\cite{j-w02, j-w03} with $n=3,$ and by Singer in \cite{si08} for the dual of  (\ref{ct1}) with $n \geqslant 0.$ For $G_n=GL(n; \mathbb F_2),$ the problem is studied by Singer~\cite{si89} for $n=2,$ and by Hung-Peterson~\cite{h-p95} for $n=3, 4.$ In case $G_n$ is the trivial group, we have 
\begin{align}\label{ct2}
\mathbb F_2 {\otimes}_{\mathcal{A}}\mathcal P_{n}=\left\{(\mathbb{F}_{2}{\otimes}_{\mathcal{A}}\mathcal P_{n})_d\right\}_{d \geqslant 0} \cong \mathcal P_{n}/{\mathcal{A}^{+}}\mathcal P_{n}, 
\end{align}
as a modular representation of $GL(n; \mathbb F_2).$ Here, $(\mathcal P_n)_d$ is the subspace of $\mathcal P_n$ consisting of all the homogeneous polynomials of degree $d$ in $\mathcal P_n,$ and $(\mathbb F_2{\otimes}_{\mathcal{A}}\mathcal P_n)_d$ is the subspace of $\mathbb F_2{\otimes}_{\mathcal{A}}\mathcal P_n$ consisting of all the classes represented by the elements in $(\mathcal P_n)_d.$

A polynomial $u$ in $\mathcal P_n$  is called {\it hit} if it can be written as a finite sum $u= \sum_{i\geqslant 0}Sq^{2^i}(f_i)$ 
for suitable polynomials $f_i$.  That means $u$ belongs to $\mathcal{A}^+\mathcal P_n$,  where $\mathcal{A}^+$ denotes the augmentation ideal in $\mathcal{A}$.

For a natural number $d,$ let $\alpha(d)$ be the number of digits 1 in the binary expansion of $d.$ We define a function $\mu : \mathbb N \longrightarrow \mathbb N$ is given by
\begin{align*} \mu (0) =0, {\text{ and }} \mu(d) &= \text{min} \{ m \in \mathbb N \ : \  d= \sum_{i = 1}^{m}(2^{d_i}-1), d_i > 0 \} \\
& =\text{min} \{ m \in \mathbb N \ : \  \alpha(d+m) \leqslant m \}.
\end{align*}

In \cite{pe87}, Peterson conjectured that as a module over the Steenrod algebra $\mathcal{A},$  $\mathcal P_n$ is generated by monomials in degree $d$ and satisfying the inequality $\alpha(d+n) \leqslant n$, where $\alpha(d)$ is the number of digits one in the binary expension of $d,$ and proved it for $n\leqslant 2,$ in general, it is proved by Wood \cite{wo89}. This is an extremely useful tool for determining $\mathcal{A}$-generators for $\mathcal P_n.$ 

One of the main tools in the study of the hit problem is Kameko's squaring operation 
$$\widetilde{Sq}^0_*:= (\widetilde{Sq}^0_*)_{(n; 2n+d)}: (\mathbb{F}_{2}{\otimes}_{\mathcal{A}} \mathcal P_n)_{2n+d} \to (\mathbb{F}_{2}{\otimes}_{\mathcal{A}}\mathcal P_n)_{d},$$
which is induced by an $\mathbb F_2$-linear map $S_n: \mathcal P_n \to \mathcal P_n$, given by
$$
S_n(x) = 
\begin{cases}y, &\text{if }x=x_1x_2\ldots x_ky^2\\  
0, & \text{otherwise} \end{cases}
$$
for any monomial $x \in \mathcal P_n$. Clearly, $ (\widetilde{Sq}^0_*)_{(n; 2n+d)}$ is an $\mathbb F_2$-epimorphism. 

From the results of Wood \cite{wo89}, Kameko \cite{ka90}, and Sum \cite{su15}, the hit problem is reduced to the case of degree $d$ of the form $d=r(2^t-1)+2^tm,$ where $r, m, t$ are non-negative intergers such that $0 \leqslant \mu(m)<r \leqslant n.$

Now, the structure of the space (\ref{ct2}) was completely calculated for $n \leqslant 4,$  (see Peterson~\cite{pe87} for $n=1,$ and $n=2,$ see Kameko for $n=3$ in his thesis \cite{ka90}, see Sum~\cite{su15} for $n = 4$). For $n \geqslant 5$, it is still unsolved, even in the case of $n=5$ with the help of computers. 

In the present paper, we develop a result in \cite{ph31} on the hit problem for $\mathcal P_n$ in the generic degree $r(2^s-1)+m.2^s,$ where $r=n=5,\ m=13,$ and $s$ is an arbitrary non-negative integer. Remarkably, for $s=0,$ and $s=1,$ this problem has been studied by Phuc \cite{ph20ta}, and \cite{ph31}, respectively. Moreover, as an application of the above results, we get the dimension results for the graded polynomial algebra in the generic degree $d=(n-1).(2^{n+u-1}-1)+\ell.2^{n+u-1}$ where $u$ is an arbitrary non-negative integer, $\ell \in \{23, 67 \},$ and $n=6.$

Note that the general linear group $GL(n; \mathbb F_2)$ acts naturally on $\mathcal P_n$ by matrix substitution. Since the two actions of $\mathcal A$ and $GL(n; \mathbb F_2)$ upon $\mathcal P_n$ commute with each other, there is an action of $GL(n; \mathbb F_2)$ on $\mathbb F_2 \otimes_{\mathcal A}\mathcal P_n.$ From this event, the hit problem becomes one of the main tools for the studying of the general linear groups over $\mathbb F_2.$ Recently, the hit problem and its applications to representations of general linear groups have been presented in the books of Walker and Wood~\cite{w-w18, w-w182}.

On the other hand, one of the major applications of hit problem is in surveying a homomorphism introduced by W. M. Singer. It is a useful tool in describing the cohomology groups of the Steenrod algebra, $Ext^{n, n+*}_{\mathcal A}(\mathbb F_2, \mathbb F_2).$ 

Let $(\mathbb F_2\otimes_{\mathcal{A}}\mathcal P_n)_{d}^{GL(n; \mathbb F_2)}$ be the subspace of $(\mathbb F_2 \otimes_{\mathcal A} \mathcal P_n)_{d}$ consisting of all the $GL(n; \mathbb F_2)$-invariant classes of degree $d,$ and $\mathbb F_2{\otimes}_{GL(n; \mathbb F_2)}PH_d((\mathbb R \mathcal P^{\infty})^n)$ be dual to $(\mathbb F_2 \otimes_{\mathcal{A}}\mathcal P_n)_{d}^{GL(n; \mathbb F_2)}.$ 

Singer \cite{si89} defined the algebraic transfer, which is a homomorphism 

$$\psi_n: \mathbb F_2{\otimes}_{GL(n; \mathbb F_2)}PH_*((\mathbb R \mathcal P^{\infty})^n) \longrightarrow Ext^{n, n+*}_{\mathcal A}(\mathbb F_2, \mathbb F_2).$$

\medskip
Singer has indicated the importance of the algebraic transfer by showing that $\psi_n$ is a isomorphism with $n=1,2$ and at some other degrees with $n=3,4$, but he also disproved this for $\psi_{5}$ at degree 9, and then gave the following conjecture.

\medskip
\begin{con} For any $n \geqslant 0,$ the algebraic transfer $\psi_n$ is a monomorphism. 
\end{con}

\medskip
It could be seen from the work of Singer the meaning and necessity of the hit problem. In \cite{bo93}, Boardman confirmed this again by using the modular representation theory of linear groups to show that $\psi_3$ is also an isomorphism. 

For $n \geqslant 4,$ the Singer algebraic transfer was studies by many authors (See Boardman \cite{bo93}, Bruner-Ha-Hung \cite{b-h-h05}, Minami \cite{mi99}, Sum-Tin \cite{s-t15} and others). However, Singer's conjecture is still open for $n \geqslant 4.$

The results of the hit problem are used to study and verify the Singer conjecture for the algebraic transfer. The behavior of the fifth Singer algebraic transfer in the degree $d_s=5(2^s -1)+13.2^s$ was also discussed at the end of this paper.

Next, in Section \ref{s2}, we recall some needed information on admissible monomials in $\mathcal P_n.$  The main results are presented in Section \ref{s3}. The proofs of the main results will be presented in Section \ref{s4}.

\section{Preliminaries}\label{s2} 
\setcounter{equation}{0}

\medskip
In this section, we recall a few useful preliminaries on the admissible monomials, spike monomials, and the hit monomials from Kameko~\cite{ka90}, and Sum \cite{su15}, which will be used in the next section.

We will denote by $\mathbb N_n = \{1,2, \ldots , n\}$ and
\begin{align*}
X_{\mathbb J} = X_{\{j_1,j_2,\ldots , j_s\}} =
 \prod_{j\in \mathbb N_n\setminus \mathbb J}x_j , \ \ \mathbb J = \{j_1,j_2,\ldots , j_s\}\subset \mathbb N_n,
\end{align*}
In particular, $X_{\mathbb N_n} =1,\
X_\emptyset = x_1x_2\ldots x_n,$ 
$X_j = x_1\ldots \hat x_j \ldots x_n, \ 1 \leqslant j \leqslant n,$ and $X:=X_n \in \mathcal P_{n-1}.$

Let $\alpha_t(d)$ be the $t$-th coefficient in dyadic expansion of $d.$ Then, $d= \sum_{t \geqslant 0} \alpha_{t}(d).2^t$ where $\alpha_t(d) \in \{0, 1\}.$
Let $x=x_1^{a_1}x_2^{a_2}\ldots x_n^{a_n} \in \mathcal P_n$. Denote $\nu_j(x) = a_j, 1 \leqslant j \leqslant n$.  
Set 
$$\mathbb J_t(x) = \{j \in \mathbb N_n :\alpha_t(\nu_j(x)) =0\},$$
for $t\geqslant 0$. Then, we have
$x = \prod_{t\geqslant 0}X_{\mathbb J_t(x)}^{2^t}.$ 

\medskip
\begin{defn}
For a monomial  $x$ in $\mathcal P_n$,  define two sequences associated with $x$ by
\begin{align*} 
\omega(x)=(\omega_1(x),\omega_2(x),\ldots , \omega_i(x), \ldots),\ \
\sigma(x) = (\nu_1(x),\nu_2(x),\ldots ,\nu_n(x)),
\end{align*} 
where
$\omega_i(x) = \sum_{1\leqslant j \leqslant n} \alpha_{i-1}(\nu_j(x))= \deg X_{\mathbb J_{i-1}(x)},\ i \geqslant 1.$
The sequence $\omega(x)$ is called  the weight vector of $x$. 
\end{defn}

\medskip
Let $\omega=(\omega_1,\omega_2,\ldots , \omega_i, \ldots)$ be a sequence of non-negative integers.  The sequence $\omega$ is called  the weight vector if $\omega_i = 0$ for $i \gg 0$.  For a  weight vector $\omega$,  we define $\deg \omega = \sum_{i > 0}2^{i-1}\omega_i$.

The sets of all the weight vectors and the exponent vectors are given the left lexicographical order. 
  
Denote by $\mathcal P_n(\omega)$ the subspace of $\mathcal P_n$ spanned by all monomials $y$ such that
$\deg y = \deg \omega$, $\omega(y) \leqslant \omega$, and by $\mathcal P_n^-(\omega)$ the subspace of $\mathcal P_n$ spanned by all monomials $y \in \mathcal P_n(\omega)$  such that $\omega(y) < \omega$. 

\medskip
\begin{defn}\label{dntd} Let $\omega$ be a weight vector and $u,\ v$ two polynomials of the same degree in $\mathcal P_n.$ We define the equivalence relations $``\equiv "$ and $``\equiv_{\omega} "$ on $\mathcal P_n$ by stating that

\medskip
${\rm (i)}$ $u \equiv v$ if and only if $u - v \in \mathcal A^+\mathcal P_n$. 

\medskip
${\rm (ii)}$ $u \equiv_{\omega} v$ if and only if $u, v \in \mathcal P_n(\omega)$ and $u-v \in \big(\mathcal A^+\mathcal P_n \cap \mathcal P_n(\omega) + \mathcal P_n^-(\omega)\big)$ denoting the set of all the $\mathcal A$-decomposable elements, in which $Sq^{2^r}: \mathcal P_{n-2^r} \to \mathcal P_n,$ and $\mathcal A^+$ is the kernel of an epimorphism $\mathcal A \to \mathbb F_2$ of graded $\mathbb F_2$-algebras. If $u$ is a linear combination of elements in the images of the Steenrod squaring operations $Sq^{2^r}$ for $r>0,$ then $u$ is called decomposable. Also, $u$ is said to be $\omega$-decomposable if $u \equiv_{\omega} 0.$
\end{defn}

\medskip
Then, we have an $\mathbb F_2$-qoutient space of $\mathcal P_n$ by the equivalence relation $``\equiv_{\omega} "$ as follows

$$Q\mathcal P_n(\omega)= \mathcal P_n(\omega)/ ((\mathcal A^+\mathcal P_n\cap \mathcal P_n(\omega))+\mathcal P_n^-(\omega)).$$ 

\medskip
\begin{defn}\label{defn3} 
Let $u,\ v$ be monomials of the same degree in $\mathcal P_n$. We say that $u <v$ if and only if one of the following holds: \ 

\medskip
${\rm (i)}$ $\omega (u) < \omega(v)$;

\medskip
${\rm (ii)}$ $\omega (u) = \omega(v),$ and and $\sigma(u) < \sigma(v).$
\end{defn}

\medskip
\begin{defn}
A monomial $u$ is said to be inadmissible if there exist monomials $v_1,v_2,\ldots, v_m$ such that $v_i<u$ for $i=1,2,\ldots , m$ and $u - \sum_{t=1}^mv_i \in \mathcal A^+\mathcal P_n.$ 
\end{defn}

\medskip
We say $u$ is admissible if it is not inadmissible. From these events, it is important to remark that the set of all the admissible monomials of degree $d$ in $\mathcal P_n$ is a minimal set of $\mathcal{A}$-generators for $\mathcal P_n$ in degree $d.$ And therefore, $(\mathbb{F}_{2}{\otimes}_{\mathcal{A}}\mathcal P_{n})_d$ is an $\mathbb F_2$-vector space with a basis consisting of all the classes represent by the elements in $(\mathcal P_n)_d.$

\medskip
\begin{defn} Let $u \in \mathcal P_n.$ We say $u$ is strictly inadmissible if and only if there exist monomials $v_1,v_2,\ldots, v_m$ such that $v_j<u,$ for $j=1,2,\ldots , m$ and $u = \sum_{j=1}^m v_j + \sum_{i=1}^{2^s-1}Sq^i(f_i)$ with $s = \max\{k : \omega_k(u) > 0\},$  $f_i \in \mathcal P_n.$
\end{defn}

\medskip
It is easy to check that if $u$ is strictly inadmissible monomial, then it is inadmissible monomial. 

\medskip
\begin{thm}[Kameko \cite{ka90}, Sum \cite{su15}]\label{dlk-s}  Suppose that $u, v, w \in \mathcal P_n$ satisfying the conditions $\omega_t(u) = 0$ if $t > k>0$, $\omega_r(w) \ne 0$ and $\omega_t(w) = 0$ if $t > r>0$. Then,  

\medskip
{\rm (i)} $uw^{2^k}$ is inadmissible if $w$ is inadmissible.

\medskip
{\rm (ii)}  $wv^{2^{r}}$ is strictly inadmissible if $w$ is strictly inadmissible.
\end{thm} 

\medskip
\begin{defn}\label{spi} Let $z=x_1^{d_1}x_2^{d_2}\ldots x_n^{d_n} \in \mathcal P_n.$ The monomial $z$ is called a spike if $d_j=2^{t_j}-1$ for $t_j$ a non-negative integer and $j=1,2, \ldots , n$. Moreover, $z$ is called the minimal spike, if it is a spike such that $t_1>t_2>\ldots >t_{r-1}\geqslant t_r>0$ and $t_j=0$ for $j>r.$
\end{defn}

\medskip
The following is a Singer's criterion on the hit monomials in $\mathcal P_n.$

\medskip
\begin{thm}[Singer~\cite{si89}]\label{dlsig} Assume that $z$ is the minimal spike of degree $d$ in $\mathcal P_n,$ and $u \in (\mathcal P_n)_d$ satisfying the condition $\mu(d) \leqslant n.$  If $\omega(u) < \omega(z),$ then $u$ is hit.
\end{thm}

\medskip
We will denote by $\mathcal P_n^0$ and $\mathcal P_n^+$ the $\mathcal{A}$-submodules of $\mathcal P_n$ spanned all the monomials $x_1^{d_1}x_2^{d_2}\ldots x_n^{d_n}$ such that $d_1 \ldots d_n=0,$ and $d_1 \ldots d_n > 0,$ respectively. It is easy to see that $\mathcal P_n^0$ and $\mathcal P_n^+$ are the $\mathcal{A}$-submodules of $\mathcal P_n.$ Then, we have a direct summand decomposition of the $\mathbb F_2$-vector spaces: 
$$\mathbb F_2{\otimes}_{\mathcal{A}}\mathcal P_n =(\mathbb F_2{\otimes}_{\mathcal{A}}\mathcal P_n^0) \oplus (\mathbb F_2{\otimes}_{\mathcal{A}}\mathcal P_n^+).$$ 

From now on, let us denote by $\mathcal C_{d}^{\otimes n}$ the set of all admissible monomials of degree $d$  in $\mathcal P_n.$ For $f \in \mathcal P_n,$ we denote by $[f]$ the class in $\mathbb F_2 {\otimes}_{\mathcal{A}}\mathcal P_{n}$ represented by $f.$ Denote by $|S|$ the cardinal of a set $S.$

\section{Statement of the main results}\label{s3} 
\setcounter{equation}{0}

\medskip
In this section, we list the main results of this paper, proofs of the main results will be presented in the next section. First, we study the hit problem for $\mathcal P_5$ in the generic degree $k_s:=5(2^s-1)+13.2^s$ with $s$ an arbitrary non-negative integer.

For $s=0,$ we have $k_0=5(2^0-1)+13.2^0=13.$ Then, $\dim(\mathbb{F}_{2}{\otimes}_{\mathcal{A}}\mathcal P_5)_{13} = 250.$ This result has been computed in \cite{ph20ta} by explicitly determining all admissible monomials of $\mathcal P_5$ in degree $13$ (see Phuc \cite{ph20ta}). 

Moreover, for $s=1,$ $k_1=5(2^1-1)+13.2^1=31,$ the following result has been shown in \cite{ph31} by Phuc.

\medskip
\begin{thm}[Phuc~\cite{ph31}]\label{dlph31} There exist exactly $866$ admissible monomials of degree $31$ in $\mathcal P_5.$ Consequently, 
$\dim (\mathbb F_2 \otimes_{\mathcal A} \mathcal P_5)_{31} = 866.$ 
\end{thm}

\medskip
One of our main contributions in this article is to study the $\mathbb F_2$-graded vector space $(\mathbb F_2 \otimes_{\mathcal A} \mathcal P_5)_{5(2^s-1)+13.2^s},$ where $s$ is any positive integer greater than one.

For $s=2,$ we have $k_2=5(2^2-1)+13.2^2=67.$ Since $\mathcal P_n = \oplus_{d\geqslant 0} (\mathcal P_n)_d$ is the graded polynomial algebra, and Kameko's homomorphism 
$$(\widetilde{Sq}^0_*)_{(5; 67)}: (\mathbb F_2 \otimes_{\mathcal A} \mathcal P_5)_{67} \longrightarrow (\mathbb F_2 \otimes_{\mathcal A} \mathcal P_5)_{31}$$
is an $\mathbb F_2$-epimorphism, it follows that
\begin{equation*}\label{ttt1} (\mathbb F_2 \otimes_{\mathcal A} \mathcal P_5)_{67} \cong (\mathbb F_2 \otimes_{\mathcal A} \mathcal P_5^0)_{67}\bigoplus \big(\text{\rm Ker}(\widetilde{Sq}^0_*)_{(5; 67)}\cap (\mathbb F_2 \otimes_{\mathcal A} \mathcal P_5^+)_{67}\big) \bigoplus \text{Im}(\widetilde{Sq}^0_*)_{(5; 67)}
\end{equation*}

Consider the homomorphism $\Phi : \mathcal P_5 \to \mathcal P_5$ is an $\mathbb F_2$-homomorphism determined by $\Phi(x)=\prod_{i=1}^5x_ix^2,$ for $x \in \mathcal P_5.$ Suppose that $\mathcal C_{31}^{\otimes 5}=\{ a_i \in \mathcal P_5 : 1 \leqslant i \leqslant 866\}.$ Then, we set 
$$\mathcal D_{Im}^{\otimes 5}=\{ [b_i] : \ b_i=\Phi(x), \text{ for all } x \in \mathcal C_{31}^{\otimes 5}\}.$$ 

\medskip
From this, we easily conclude the following theorem.

\medskip
\begin{thm}\label{dl1} $\text{Im}(\widetilde{Sq}^0_*)_{(5; 67)}$ is isomorphic to a subspace of $(\mathbb F_2 \otimes_{\mathcal A} \mathcal P_5)_{67}$ generated by all the classes of the form $[b_i]$, with $1 \leqslant i \leqslant 866.$ That means, the set $\mathcal D_{Im}^{\otimes 5}$ is a basis of the $\mathbb F_2$-vector space $\text{Im}(\widetilde{Sq}^0_*)_{(5; 67)}.$ 
\end{thm}

\medskip
For any $1 \leqslant t \leqslant 5$, consider the homomorphism $\mathcal T_t: \mathcal P_{4} \to \mathcal P_5$ of algebras by substituting
$$\mathcal T_t(x_i) = \begin{cases} x_i, &\text{ if } 1 \leqslant i <t,\\
x_{i+1}, &\text{ if } t \leqslant i <5.
\end{cases}$$
Then, $\mathcal T_t$ is a homomorphism of $\mathcal A$-modules. Put $\mathcal D_{0}^{\otimes 5}:=\{ c_i : \ c_i \in \bigcup _{t=1}^5\mathcal T_t(\mathcal C_{67}^{\otimes 4})\}.$

Using the result in Sum~\cite{su15} (see Theorem 1.4), an easy computation shows that $|\mathcal D_{0}^{\otimes 5}|=460.$ More specifically, we have the following.

\medskip
\begin{thm}\label{dl2} The set $\{ [c_i] : \ 1\leqslant i \leqslant 460 \}$ is a basis of the $\mathbb F_2$-vector space $(\mathbb F_2 \otimes_{\mathcal A} \mathcal P_5^0)_{5(2^2-1)+13.2^2}.$ Consequently, $\dim(\mathbb F_2 \otimes_{\mathcal A} \mathcal P_5^0)_{5(2^2-1)+13.2^2}=460.$
\end{thm} 

\medskip
\begin{rem}\label{nx1} We recall a result in Mothebe-Kaelo-Ramatebele \cite{m-k-o16} as follows: 

\medskip
Set $\mathcal N_{(n,t)}=\{I= (i_1, i_2, \ldots, i_t) : 1 \leqslant i_1 <\ldots < i_t \leqslant n\}$, $1 \leqslant t < n$. For $I\in \mathcal N_{(n,t)}$, define the homomorphism $f_I : \mathcal P_t \to \mathcal P_n$ of algebras by substituting $f_I(x_{\ell}) = x_{i_{\ell}}$ with $1 \leqslant \ell \leqslant t$. Then, $f_I$ is a monomorphism of $\mathcal A$-modules. We have a direct summand decomposition of the $\mathbb F_2$-vector subspaces: 
\medskip
$$ \mathbb F_2 \otimes_{\mathcal A} \mathcal P_n^0  = \bigoplus\limits_{1\leqslant t\leqslant n-1}\bigoplus\limits_{I \in \mathcal N_{(n,t)}}(Qf_I(\mathcal P_t^+)),$$
where $Qf_I(\mathcal P_t^+) = \mathbb{F}_2\otimes_{\mathcal{A}}f_I(\mathcal P_t^+)$. 

\medskip
Hence, $\dim (Qf_I(\mathcal P_t^+))_d = \dim (\mathbb F_2 \otimes_{\mathcal A} \mathcal P_t^+)_d$ and $|\mathcal N_{(n,t)}| = \binom n t$. Combining with the results in Wood \cite{wo89}, one gets 

\medskip
\begin{equation*}\label{ctscbs} \dim ( \mathbb F_2 \otimes_{\mathcal A} \mathcal P_n^0)_d = \sum\limits_{\mu(d) \leqslant t \leqslant n-1}\binom{n}{t}\dim(\mathbb F_2 \otimes_{\mathcal A} \mathcal P_t^+)_d.
\end{equation*} 

\medskip
Since $\mu(67)=3,$ it folows that if $t < 3$ then the spaces $(\mathbb F_2 \otimes_{\mathcal A} \mathcal P_t^+)_{67}$ are trivial. Moreover, by using the results in Kameko~\cite{ka90} and Sum~\cite{su15}, we have 
$$ \dim(\mathbb F_2 \otimes_{\mathcal A} \mathcal P_t^+)_{67} = 
\begin{cases}14, &\text{if } \ t=3,\\  
64, & \text{if } \ t=4.
\end{cases}
$$

From this, we get $$\dim (\mathbb F_2 \otimes_{\mathcal A} \mathcal P_5^0)_{67} = \binom{5}{3}\dim(\mathbb F_2 \otimes_{\mathcal A} \mathcal P_3^+)_{67} +\binom{5}{4}\dim(\mathbb F_2 \otimes_{\mathcal A} \mathcal P_4^+)_{67}=460.$$
\end{rem}

Next, we explicitly determine the vector space $\text{\rm Ker}(\widetilde{Sq}^0_*)_{(5; 67)}\cap (\mathbb F_2 \otimes_{\mathcal A}\mathcal P_5^+)_{67}.$ We will denote by $Q\mathcal P_n^+(\omega)=Q\mathcal P_n(\omega) \cap (\mathbb F_2{\otimes}_{\mathcal{A}}\mathcal P_n^+).$ Then, we have the following theorem.

\medskip
\begin{thm}\label{dl3} Let $\widetilde{\omega_1}:=(3, 4, 2, 2, 2),$\ $\widetilde{\omega_2}:=(3, 4, 4, 3, 1),$ and $ \widetilde{\omega_3}:=(3, 2, 1, 1, 1).$ 

\medskip
{\rm (i)} Suppose $u$ is an element of $(\mathcal C_{67}^{\otimes 5}\cap \mathcal P_5^+)$ such that $[u]$ does not belong to $\text{\rm Im}(\widetilde{Sq}^0_*)_{(5; 67)},$ then $\omega(u)$ is one of  the following sequences: $\widetilde{\omega_m}$, with $m=1, 2, 3.$
Moreover, we have an isomorphism of the $\mathbb F_2$-vector spaces: 
$$\text{\rm Ker}(\widetilde{Sq}^0_*)_{(5; 67)}\cap (\mathbb F_2{\otimes}_{\mathcal{A}}\mathcal P_5^+)_{67} \cong \bigoplus_{m=1}^{3} Q\mathcal P_5^+(\widetilde{\omega_m}).$$  

{\rm (ii)} We have $\dim \text{\rm Ker}(\widetilde{Sq}^0_*)_{(5; 67)}\cap (\mathbb F_2{\otimes}_{\mathcal{A}}\mathcal P_5^+)_{67}=161.$
\end{thm}

\medskip
From the results of Theorems ~\ref{dl1}, \ref{dl2}, and~ \ref{dl3}, we obtain the following corollary.

\medskip
\begin{corl}\label{corl1} There exist exactly $1487$ admissible monomials of degree $67$ in $\mathcal P_5.$ Consequently, 
$\dim (\mathbb F_2 \otimes_{\mathcal A} \mathcal P_5)_{5(2^2-1)+13.2^2} = 1487.$ 
\end{corl}

\medskip
Consider the degree $k_{s}=5(2^{s} -1)+13.2^{s},$ for any $s \geqslant 3.$ We set 

$$\xi(n; d) = \max\{0,n- \alpha(d+n) -\zeta(d+n)\},$$ 

\medskip
\noindent{where $\zeta(u)$ the greatest integer $v$ such that $u$ is divisible by $2^v.$}
We recall a result in \cite{t-s16} the following. 

\medskip
\begin{thm}[Tin-Sum~\cite{t-s16}]\label{dlt-s16}Let $d$ be an arbitrary non-negative integer. Then

$$(\widetilde{Sq}^0_*)^{r-t}: (\mathbb F_2 \otimes_{\mathcal A} \mathcal P_n)_{n(2^r-1) + 2^rd} \longrightarrow (\mathbb F_2 \otimes_{\mathcal A} \mathcal P_n)_{n(2^t-1) + 2^td}$$

\medskip
\noindent is an isomorphism of $GL(n; \mathbb F_2)$-modules for every $r \geqslant t$ if and only if $t \geqslant  \xi(n; d).$
\end{thm}

\medskip
As an application of Theorem \ref{dlt-s16}, we get the following theorem.

\medskip
\begin{thm}\label{dl4} The set $\big\{ [x]: x \in \Phi^{s-2}\big( \mathcal C_{5(2^2-1) + 13.2^2}^{\otimes 5}\big) \big\}$ is a basis of the $\mathbb F_2$-vector space $(\mathbb F_2 {\otimes}_{\mathcal{A}}\mathcal P_5)_{5(2^{s} -1)+13.2^{s}}$, for any $s > 2.$ This implies that the $\mathbb F_2$-vector space
$(\mathbb F_2 \otimes_{\mathcal A} \mathcal P_5)_{5(2^{s} -1)+13.2^{s}}$ has dimension $1487,$ for any $s > 2.$
\end{thm}

\medskip
Now, we describe the next main results by studying the $\mathbb F_2$-graded vector space $\mathbb F_2 {\otimes}_{\mathcal{A}}\mathcal P_n$ in some generic degrees of the form ${(n-1)(2^{r+4} -1)+m.2^{r+4}},$ where $r$ is an arbitrary positive integer, $m \in \{23, 67\},$ and in the case $n=6.$

As is well known, after explicitly determining $\mathbb F_2 \otimes_{\mathcal A} \mathcal P_4,$ Sum \cite{su15} has established an inductive formula by $n$ for the dimension of the vector space $(\mathbb F_2 \otimes_{\mathcal A} \mathcal P_n)_{d},$ where $d$ is of general degree (see Theorem 1.3).

From the results in Sum~\cite{su15}, combined with the above results we get the following theorem.

\medskip
\begin{thm}\label{dl5} The vector space $(\mathbb F_2 \otimes_{\mathcal A} \mathcal P_6)_{5(2^{s+4}-1)+67.2^{s+4}}$ is $93681$-dimensional, for all positive integers $s.$
\end{thm}

On the other hand, based on the results in \cite{ti21aejm} for the hit problem of five variables, we also obtain the following.

\begin{thm}\label{dl6} For any integer $r > 0$, there exist exactly $78435$ admissible monomials of degree $5.2^{r+4} - 23.2^{r+4}$ in $\mathcal P_6$. Consequently, $\dim (\mathbb F_2 \otimes_{\mathcal{A}} \mathcal P_6)_{5.2^{r+4} - 23.2^{r+4}} = 78435,$ for all $r>0.$
\end{thm}

\medskip
\section{Proofs of main results}\label{s4}
\setcounter{equation}{0}
\medskip

\subsection{\rm{Proof of Theorem \ref{dl3}}}\label{s4.1} \

\medskip

We first prove Part {\rm (i)} of the theorem. Let us denote by 

$$Q\mathcal P_5^\omega := \text{ Span }\{[x] \in \mathbb F_2{\otimes}_{\mathcal A}\mathcal P_5 : x \text{ \rm  is admissible and } \omega(x) = \omega\}.$$  

\medskip

Using the results in Walker-Wood~\cite{w-w18}, we see that

\begin{equation}\label{ww1}
(\mathbb F_2{\otimes}_{\mathcal A}\mathcal P_5)_d = \bigoplus_{\deg \omega = d}Q \mathcal P_5^\omega \cong \bigoplus_{\deg \omega = d}Q \mathcal P_5(\omega),
\end{equation}

Since $\mathcal P_n= \oplus_{d\geqslant 0} (\mathcal P_n)_d$ is the graded polynomial algebra, combined with the above results, it follows that 
$(\mathbb F_2{\otimes}_{\mathcal A}\mathcal P_5^+)_{67} = \bigoplus_{\deg \omega = 67}Q \mathcal P_5^+(\omega).$

Suppose that $u$ is an admissible monomial of degree $67$ in $\mathcal P_5^+$ such that $[u]$ belongs to $\text{\rm Ker}(\widetilde{Sq}^0_*)_{(5; 67)}.$ 
It is easy to check that $z = x_1^{63}x_2^{3}x_3$ is the minimal spike of degree $67$ in $\mathcal P_5$ and  $\omega(z) = \widetilde{\omega_3}.$ Using Theorem \ref{dlsig}, we obtain $\omega_1(u) \geqslant \omega_1(z) =3$. Since the degree of $u$ is odd, it follows that either $\omega_1(u) =3$ or $\omega_1(u) =5.$ 

If $\omega_1(u) =5,$ then $u=\prod_{i=1}^{5}x_iv^2$ with $v$ a monomial in $(\mathcal P_5)_{31}.$ By Theorem \ref{dlk-s}, $v$ is an admissible monomial and $[v] \not = 0$. Thus, $[v]=\text{\rm Ker}(\widetilde{Sq}^0_*)_{(5; 67)}([u]) \not = 0.$ This contradicts the fact that $[u]$ belongs to $\text{\rm Ker}(\widetilde{Sq}^0_*)_{(5; 67)}.$ From this, one gets $\omega_1(u) =3.$ Then, we have $u=x_ix_jx_\ell y^2$ with $1 \leqslant i < j < \ell \leqslant 5,$ where $y \in (\mathcal P_5)_{32}.$  Using Theorem \ref{dlk-s}, $y$ is also admissible. By a simple computation shows that $\omega(y) = (4, 4, 3, 1)$ or $ \omega(y) =  (4, 2, 2, 2)$ or $\omega(y) =  (2, 1, 1, 1, 1).$
And therefore, we get $\omega(u)=\widetilde{\omega_m}$, with $m=1, 2, 3.$

From these above, combined with the aid of (\ref{ww1}), we have a direct summand decomposition of the $\mathbb F_2$-vector spaces:
$$\text{\rm Ker}(\widetilde{Sq}^0_*)_{(5; 67)}\cap (\mathbb F_2{\otimes}_{\mathcal A}\mathcal P_5^+)_{67} =\bigoplus_{m=1}^{3} Q\mathcal P_5^+(\widetilde{\omega_m}).$$ Part {\rm (i)} is proved.

The proof of Part {\rm (ii)} of the above theorem is too long and computationally very technical. We sketch its proof as follows: 

Consider the set $B_{67}^{\otimes 5}:=\{ X_{\{i, j\}}.F^2 : 1 \leqslant i<j \leqslant 5, F \in \mathcal C_{32}^{\otimes 5} \}.$ By using Theorem \ref{dlk-s}, we see that if $X$ is an admissible monomial of degree $67$ in $\mathcal P_5$, then $X \in B_{67}^{\otimes 5}.$  Putting $$\mathcal C_{67}^{\otimes 5}(\omega):=\{ X \in B_{67}^{\otimes 5}: X \text{ admissible monomial, and } \omega(X)=\widetilde{\omega_m},  \text{ for }\ m=1, 2, 3 \}.$$

By direct calculations, we filter out and eliminate the inadmissible monomials in $B_{67}^{\otimes 5},$ one gets $|\mathcal C_{67}^{\otimes 5}(\omega)|=161.$ 

Note that the result dimension of $\text{\rm Ker}(\widetilde{Sq}^0_*)_{(5; 67)}\cap (\mathbb F_2{\otimes}_{\mathcal A}\mathcal P_5^+)_{67}$ has been verified by using a computer calculation program in SAGE by H-V Vu. We would like to say thank you for his support.

\subsection{\rm{Proof of Theorem \ref{dl4}}}\label{s4.2} \

\medskip

It is easy to check that for $n=5$ and $d=67$ then $\alpha(d+n)=\alpha(72) = 3,$ and $\zeta(d+n)=\zeta(72)=3,$ and therefore $\xi(n; d)=0.$ By Theorem \ref{dlt-s16}, we get an isomorphism of $\mathbb F_2$-graded vector spaces:
$$(\mathbb F_2 \otimes_{\mathcal A} \mathcal P_5)_{5(2^r-1) + 2^r.67} \cong (\mathbb F_2 \otimes_{\mathcal A} \mathcal P_5)_{5(2^0-1) + 2^0.67}$$ for all $r\geqslant 0.$ 

And therefore, we obtain $$(\mathbb F_2 \otimes_{\mathcal A} \mathcal P_5)_{5(2^s-1) + 13.2^s} \cong (\mathbb F_2 \otimes_{\mathcal A} \mathcal P_5)_{5(2^2-1) + 13.2^2}\ \text{for all } s > 2.$$ 
That means,\ $\dim (\mathbb F_2 \otimes_{\mathcal A} \mathcal P_5)_{5(2^s-1) + 13.2^s}= \dim (\mathbb F_2 \otimes_{\mathcal A} \mathcal P_5)_{5(2^2-1) + 13.2^2}=1487,$\ for all $s > 2.$ 

Moreover, for every $s > 2,$ the set $\big\{ [x]: x \in \Phi^{s-2}\big( \mathcal C_{5(2^2-1) + 13.2^2}^{\otimes 5}\big) \big\}$ is a basis of the $\mathbb F_2$-vector space $(\mathbb F_2 {\otimes}_{\mathcal{A}}\mathcal P_5)_{5(2^{s} -1)+13.2^{s}}.$

\subsection{\rm{Proof of Theorems \ref{dl5} and \ref{dl6}}}\label{s4.3} \

\medskip

Consider the degree $d= (n-1)(2^r-1) + 2^rm,$ with $r, m$ positive integers such that $1 \leqslant n-3 \leqslant \mu(m) \leqslant n-2.$ If $r \geqslant n-1,$ then
 we have the following result, which was shown in Sum~\cite{su15}.

$$ \dim(\mathbb F_2{\otimes}_{\mathcal A}\mathcal P_n)_d =(2^n-1) \dim(\mathbb F_2{\otimes}_{\mathcal A}\mathcal P_{n-1})_m.$$

\medskip
For $n=6,$ and $m=67,$ we can easily see that $\mu(67)= 3= \alpha(67+\mu(67))= \alpha(70).$ Hence, using the above result, one gets 

$$\big|\mathcal C_{5(2^{r}-1)+67.2^{r}}^{\otimes 6}\big|=63.\big|\mathcal C_{67}^{\otimes 5}\big|=93681,\   \text{for any integer}\  r \geqslant n-1=5.$$ 

\medskip
So, there exist exactly $93681$ admissible monomials of degree $5(2^{r}-1)+67.2^{r}$ in $\mathcal P_6,$ for any integer $r >4.$ Consequently, 
$\dim (\mathbb F_2 \otimes_{\mathcal A} \mathcal P_6)_{5(2^{s+4}-1)+67.2^{s+4}} = 93681,$ for all positive integers $s.$ Theorem \ref{dl5} is proved.

Moreover, Tin~\cite{ti21aejm} showed that $(\mathbb F_2{\otimes}_{\mathcal{A}}\mathcal P_5)_{5(2^{r} -1)+9.2^{r}}$ is an $\mathbb F_2$-vector space of dimension 1245, for any integer $r>0$. Consequently, $\big|\mathcal C_{5(2^{r}-1)+9.2^{r}}^{\otimes 5}\big|=1245,$ for any integer $r>0.$

Consider the degree $5(2^{u} -1)+23.2^{u}.$ It is easy to check that 

$$\mu(23)=\alpha(23 + \mu(23))=\alpha(26)=n - 3=3.$$ 

\medskip
And therefore, using the result in Sum~\cite{su15} (see Theorem 1.3), we get 

$$\big|\mathcal C_{5(2^{u}-1)+23.2^{u}}^{\otimes 6}\big|=(2^6-1)\big|\mathcal C_{23}^{\otimes 5}\big|, \   \text{for any integer}\  u \geqslant 5.$$ 

\medskip
This implies that the $\mathbb F_2$-graded vector space
$(\mathbb F_2 \otimes_{\mathcal A} \mathcal P_6)_{5(2^{v+4}-1)+23.2^{v+4}}$ has dimension $78435,$ for any integer $v >0.$ Theorem \ref{dl6} is proved.

\medskip

\section{On behavior of the Singer algebraic transfer}\label{s5}
\setcounter{equation}{0}
\medskip

We first recall the definition of the Singer algebraic transfer, which is a homomorphism  

$$\psi_n: \mathbb F_2{\otimes}_{GL(n; \mathbb F_2)}PH_*((\mathbb R \mathcal P^{\infty})^n) \longrightarrow Ext^{n, n+*}_{\mathcal A}(\mathbb F_2, \mathbb F_2).$$ 

\medskip
Here, $(\mathbb F_2\otimes_{\mathcal{A}}\mathcal P_n)_{d}^{GL(n; \mathbb F_2)}$ is the subspace of $(\mathbb F_2 \otimes_{\mathcal A} \mathcal P_n)_{d}$ consisting of all the $GL(n; \mathbb F_2)$-invariant classes of degree $d,$ and $\mathbb F_2{\otimes}_{GL(n; \mathbb F_2)}PH_d((\mathbb R \mathcal P^{\infty})^n)$ be dual to $(\mathbb F_2 \otimes_{\mathcal{A}}\mathcal P_n)_{d}^{GL(n; \mathbb F_2)}.$ 

Noting that in the dual case, we also have an algebraic homomorphism called Singer's algebraic transfer, $Tr_n:=(\psi_n)^*,$ which is a homomorphism from the homology of the Steenrod algebra, $\text{Tor}^{\mathcal A}_{n, n+d} (\mathbb F_2,\mathbb F_2),$ to the subspace of $(\mathbb{F}_{2}{\otimes}_{\mathcal{A}}\mathcal P_{n})_d$ consisting of all the $GL(n; \mathbb F_2)$-invariant classes (see Singer~\cite{si89}).

In \cite{s-t15}, and ~\cite{ti21aejm}, we based on the the results for the hit problem to study and verify the Singer conjecture for the algebraic transfer in degrees $5(2^s-1)+2^sm,$ for $s$ an arbitrary positive integer, and $m \in \{ 1,2,3 \}$, we obtain the following theorem. 

\medskip
\begin{thm}\label{dl7} {\it Let $s$ be an arbitrary positive integer. Singer's conjecture is true for $n=5$ and the generic degrees $d_s=5(2^s-1)+2^sm,$ where $m \in \{1,2,3\}.$}
\end{thm} 

\medskip
In the current article, by using Theorem \ref{dlt-s16}, we also see that $$(\mathbb F_2 \otimes_{\mathcal A} \mathcal P_5)_{5(2^{s} -1)+13.2^{s}}^{GL(5; \mathbb F_2)} \cong (\mathbb F_2  \otimes_{\mathcal A} \mathcal P_5)_{5(2^{2} -1)+13.2^{2}}^{GL(5; \mathbb F_2)},\ \text{ for all } s > 1.$$

And therefore, we get 
\medskip
$$\mathbb F_2{\otimes}_{GL(5; \mathbb F_2)}PH_{5(2^{s} -1)+13.2^{s}}((\mathbb R \mathcal P^{\infty})^5) \cong \mathbb F_2{\otimes}_{GL(5; \mathbb F_2)}PH_{5(2^{2} -1)+13.2^{2}}((\mathbb R \mathcal P^{\infty})^5),$$ 

\noindent for all $s > 1.$ Hence, we need only to compute the dimension of vector space $\mathbb F_2{\otimes}_{GL(5; \mathbb F_2)}PH_{5(2^{s} -1)+13.2^{s}}((\mathbb R \mathcal P^{\infty})^5)$ for $s=2.$ 

Noting that for $s=1,$ then $\dim\big((\mathbb F_2 \otimes_{\mathcal A} \mathcal P_5)_{5(2^{1} -1)+13.2^{1}}^{GL(5; \mathbb F_2)}\big) =2$ (see Phuc~\cite{ph31}). Moreover, since Kameko's homomorphism 
$$(\widetilde{Sq}^0_*)_{(5; 67)}: (\mathbb F_2 \otimes_{\mathcal A} \mathcal P_5)_{5(2^{2} -1)+13.2^{2}} \longrightarrow (\mathbb F_2 \otimes_{\mathcal A} \mathcal P_5)_{5(2^{1} -1)+13.2^{1}}$$

\medskip
\noindent is also an $GL(5; \mathbb F_2)$-epimorphism, it follows that
\medskip
$$\dim\big(\mathbb F_2{\otimes}_{GL(5; \mathbb F_2)}PH_{5(2^{s} -1)+13.2^{s}}((\mathbb R \mathcal P^{\infty})^5)\big) \leqslant \dim(\widetilde{Sq}^0_*)^{GL(5; \mathbb F_2)}_{(5; 67)} +2,  \text{ for all } s > 1.$$

\medskip
Remarkably, for any $s>1,$ all elements of  $\mathbb F_2{\otimes}_{GL(5; \mathbb F_2)}PH_{5(2^{s} -1)+13.2^{s}}((\mathbb R \mathcal P^{\infty})^5)$ are of the form $\big(\xi[\Phi(\upsilon)]+\gamma[\Phi(\beta)]+[f]\big)^*,$ where $\xi, \gamma \in \mathbb F_2,$ $\upsilon$ and $\beta$ are two generators of $(\mathbb F_2 \otimes_{\mathcal A} \mathcal P_5)_{5(2^{1} -1)+13.2^{1}}^{GL(5; \mathbb F_2)},$ and $f \in (\mathcal P_5)_{5(2^{2} -1)+13.2^{2}}$ such that $[f]$ belongs to $\text{\rm Ker}(\widetilde{Sq}^0_*)_{(5; 67)}.$ By using techniques of the hit problem, combining the computation of $\text{\rm Ext}_{\mathcal A}^{n, n+*}(\mathbb Z_2, \mathbb Z_2)$ by Lin \cite{li08} and Chen \cite{ch11}, we will describle explicitly all these elments, and verify the Singer's conjecture for the fifth algebraic transfer in the near future.

\medskip
\noindent
{\bf Acknowledgment.} I would like to express my warmest thanks to my adviser, Prof. Nguyen Sum (Sai Gon University) for helpful conversations. I also thank Dr. Dang Vo Phuc (Khanh Hoa University) for a helpful discussion. This research is supported by Ho Chi Minh City University of Technology and Education (HCMUTE), Vietnam.

{}

\vskip0.7cm

\end{document}